\documentstyle[12pt]{article}
\begin{document}
\title{Tangent scrolls in prime Fano threefolds}
\author{Atanas Iliev \and Carmen Schuhmann}
\date{}
\maketitle
\begin{abstract}
{\footnotesize
In this paper we prove that any smooth prime Fano threefold,
different from the Mukai-Umemura threefold $X_{22}'$,
contains a 1-dimensional family of intersecting lines.
Combined with a result in [Sch] this implies that any morphism
from a smooth Fano threefold of index 2 to a smooth Fano threefold
of index 1 must be constant, which gives an answer in dimension 3
to a question stated by Peternell.
}
\end{abstract}
\centerline{\bf \S \ 1. \ Introduction}

\medskip 

{\bf (1.1)} \
A smooth projective variety $X$ is called a Fano variety
if the anticanonical bundle $-K_X$ is ample. Then the
index of $X$ is the largest positive integer $r = r(X)$
such that $-K_X = rH$ for some line bundle $H$ on $X$.

The smooth Fano threefold $X = X_d \subset {\bf P}^{g+1}$
($d = deg \ X$)
is called {\it prime} if ${\rho}(X) = rank \ {\bf Pic}(X) = 1$,
$r(X) = 1$, and $-K_X$ is the hyperplane bundle on $X$.
By the classification of Fano threefolds
smooth prime Fano threefolds exist {\it iff}
$3 \le g \le 12 (g \not= 11)$, and then $d = 2g-2$
(see [I1]).

\medskip 

{\bf (1.2)}
(see \S 4.2, \S 4.4 in [IP], or \S 1 in [I2]). \
Let $l$ be a line on the smooth prime Fano threefold $X$,
and let $N_{l/X}$ be the normal bundle of $l \subset X$.
Then

{\bf (1).} \
{\it either}
{\bf (a).} $N_{l/X} = {\cal O} \oplus {\cal O}(-1)$;
{\it or}
{\bf (b).} $N_{l/X} = {\cal O}(1) \oplus {\cal O}(-2)$.

{\bf (2).} \
The Hilbert scheme ${\cal H}_X$ of lines on $X$
is non-empty, any irreducible component
${\cal H}_o$ of ${\cal H}_X$ is one-dimensional,
and
{\it either}
${\cal H}_o$ is {\it non-exotic}, i.e.
$N_{l/X}$ is of type {\bf (1)(a)} for the general $l \in {\cal H}_o$;
\ {\it or}
${\cal H}_o$ is {\it exotic}, i.e.
$N_{l/X}$ is of type {\bf (1)(b)} for any $l \in {\cal H}_o$.

{\bf (3).} \
The component ${\cal H}_o$ is exotic if
{\it either}
the elements $l \in {\cal H}_o$ sweep out the tangent scroll
$R_o \subset X$ to an irreducible curve $C \subset X$;
\ {\it or}
$g = 3$ (i.e. $X = X_4$ is a quartic threefold), and then the lines
$l \in {\cal H}_o$ sweep out a hyperplane section $R_o \subset X_4$
which is a cone over a plane quartic curve, centered at some point
$x \in X_4$.

\medskip 

{\bf (1.3)} \
For example, the scheme ${\cal H}_X$ of the Fermat quartic
$X = X_4 = (x_0^4 +...+ x_4^4 = 0)$, which is a prime Fano threefold
of $g = 3$, is a union of $40$ double components each of which is
of type {\bf (1.2)(1)(b)} (see Rem. (3.5)(ii) in [I1]).

The only known example of a prime Fano threefold $X$
of $g \ge 4$ such that ${\cal H}_X$ has an exotic component
${\cal H}_o$, is the Mukai-Umemura threefold $X_{22}'$.
The scheme ${\cal H}_{X'_{22}} = 2{\cal H}_o$
and the surface $R_o$ is the hyperplane section
of $X_{22}'$ swept out by the tangent lines to a rational normal
curve $C_{12} \subset X_{22}'$ of degree $12$ (see [MU]).

\medskip 

{\bf (1.4)} \
By a theorem of Kobayashi and Ochiai
the index $r = r(Y)$ of a smooth Fano $n$-fold $Y$ can't be
greater than $n+1$; and the only smooth Fano $n$-folds
of $r \ge n$ are ${\bf P}^n$ for which $r = n+1$
and the $n$-dimensional quadric $Q^n_2$ for which $r = n$
(see e.g. [Pe],p. 106).
In particular, except ${\bf P}^3$ and $Q^3_2$, any
smooth Fano $3$-fold must have index $r \le 2$.

It is shown by Remmert and Van de Ven (for $n = 2$) and later by
Lazarsfeld (for any $n$) that the projective space ${\bf P}^n$
does not admit surjective morphisms to a smooth projective
$n$-fold $X \not= {\bf P}^n$ (see [RV], [L]).
The same is true for morphisms
$f: Q^n_2 \rightarrow X \not= {\bf P}^n, Q^n_2$ (see [PS]).
In particular, ${\bf P}^3$ and $Q^3_2$ do not admit surjective
moprphisms to smooth Fano threefolds $X$
of {\it smaller} index $r(X)$.

Let $f:Y \rightarrow X$ be a non-constant morphism between
smooth Fano $3$-folds of $\rho =  1$. By Kor. 1.5 in [RV],
${\rho}(Y)=1$ implies that $f$ must be surjective, and by
the preceding $r(Y)$ can't be $\ge 3$. Therefore
$r(Y) = 2$, $r(X) = 1$.
This gives rise to the following question stated
originally by Peternell (see (2.12)(2) in [Pe]).

\medskip

{\bf (Pe)}
{\bf Question}. \ 
{\sl
Are there non-constant (hence surjective)
morphisms $f:Y \rightarrow X$
from a smooth Fano $3$-fold $Y$ of ${\rho}(Y) = 1$
and $r(Y) = 2$ to a smooth Fano $3$-fold $X$
of ${\rho}(X) = 1$ and $r(X) = 1$ ?
}

\medskip

In this paper we give the expected negative answer to {\bf (Pe)}.

\medskip

Let $f:Y \rightarrow X$ be as above,
and assume that $f$ is non-constant.
Then $f$ must be surjective and finite since
${\rho}(Y) = 1$ (see Kor. 1.5 in [RV]). Therefore
$f^*:H^3(X,{\bf C}) \rightarrow H^3(Y,{\bf C})$
will be an embedding, in particular $h^3(X) \le h^3(Y)$
(see also [Sch]).
For any Fano threefold $h^{3,0} = 0$ and $h^3 = 2h^{2,1}$
since the anticanonical class is ample.
Therefore $h^{2,1}(X) \le h^{2,1}(Y)$.
Since $h^{2,1}(Y) \le 21$ for any Fano 3-fold $Y$ of $r = 2$
(see [I1]), then the existence of a non-constant morphism
$f:X \rightarrow Y$ as in {\bf (Pe)} implies that
$h^{2,1}(X) \le 21$. This gives a negative answer to {\bf (Pe)}
whenever $h^{2,1}(X) > 21$.

The only smooth non-prime Fano threefolds of ${\rho} = 1$ and $r = 1$
are the sextic double solid $X_2'$ and the double quadric $X_4'$
for which the answer to {\bf (Pe)} is negative since
$h^{2,1}(X_2') = 52 > 21$ and $h^{2,1}(X_4') = 30 > 21$.
By the same argument the answer to {\bf (Pe)} is
negative also for the quartic threefold $X_4$
since $h^{2,1}(X_4) = 30 > 21$.
Any other smooth prime Fano threefold
$X = X_{2g-2} \subset {\bf P}^{g+1}, 4 \le g \le 12, g \not= 11$
has $h^{2,1}(X) \le 20$ (see [I1]).

In [Sch] is given a negative answer to {\bf (Pe)} provided
$X$ contains a conic of rank 2 (a pair of intersecting lines).
The only known Fano threefold $X$ of ${\rho}(X) = 1$
and $r(X) = 1$ without intersecting lines is the Mukai-Umemura
threefold $X = X_{22}'$, and a negative answer to {\bf (Pe)}
in case $X = X_{22}'$ is given by E. Amerik (see [Sch]).

\medskip

Therefore, in order to give a negative answer
to {\bf (Pe)}, it is enough to prove the following

\medskip

{\bf (B) Proposition.}
{\sl
Any smooth prime Fano threefold $X_{2g-2} \subset {\bf P}^{g+1}$,
$4 \le g \le 12, g \not= 11$, different from the Mukai-Umemura
threefold $X_{22}'$, contains a 1-dimensional family of conics of
rank $2$.
}

\medskip

In Section 2 we prove Proposition {\bf (B)}
for $4 \le g \le 9$ on the base of the following
technical

\medskip

{\bf (A) Lemma.}
{\sl
A smooth prime Fano threefold
$X = X_{2g-2} \subset {\bf P}^{g+1}, 3 \le g \le 9$
can't contain the tangent scroll $S_{2g-2}$
to a rational normal curve $C_g$ of degree $g$.
}

\medskip

By a result of Yu. Prokhorov,
the only smooth prime $X = X_{2g-2}$, $g = 10,12$ such that
the scheme of lines ${\cal H}_X$ has an exotic component is
the Mukai-Umemura threefold $X_{22}'$ (see [Pr]).
This implies Proposition {\bf (B)} for $g = 10,12$.
Indeed, let $X = X_{2g-2}$ be a smooth prime Fano threefold
such that the scheme of lines ${\cal H}_X$ on $X$ has
a non-exotic component ${\cal H}_o$.
Then, by Lemma (3.7) in [I1], the general element of
${\cal H}_o$ will represent a line $l \subset X$
which intersects at least one other line on $X$.

This completes the proof of Proposition {\bf (B)},
which yields a negative answer to {\bf (Pe)}.

\medskip

In Section 3 we prove Lemma {\bf (A)} for any
particular value of $g$, $3 \le g \le 9$.

For the prime Fano threefolds
$X_{2g-2} \subset {\bf P}^{g+1}$ ($g = 3,4,5,6,8$)
we prove Lemma {\bf (A)}
by using the Mukai's representation (see {\bf (3.1)})
of the smooth $X_{2g-2}$, $3 \le g \le 10$ as a complete
intersection in a homogeneous or almost-homogeneous
variety ${\Sigma}(g)$.
More concretely we see that if the
threefold $X \subset {\Sigma}(g)$ ($g = 3,4,5,6,8$)
is a complete intersection in ${\Sigma}(g)$
of the same type as the smooth prime $X_{2g-2}$,
and if $X$ contains the tangent scroll $S_{2g-2}$
to the rational normal curve $C_g$ of degree $g$,
then $X$ must be singular
-- see {\bf (3.4)}, {\bf (3.5)}, {\bf (3.6)}, {\bf (3.13)},
{\bf (3.15)}, {\bf (3.19)}.

For $g = 7$ we use the properties of the projection
from a special line $l \subset X_{2g-2}$, $g \ge 7$
to reduce the proof of Lemma {\bf (A)} for $g=7$
to the already proved Lemma {\bf (A)} for $g = 5$
-- see {\bf (3.20) - (3.24)}.
To prove Lemma {\bf (A)} in case $g = 9$ 
we can use the same approach as for $g = 7$.
But a more elegant proof, based on the
description of the double projection from a line,
had been suggested by the Referee -- see {\bf (3.25)}.

\bigskip

\bigskip 

\centerline{\bf \S \ 2. \ Lemma (A) $\Rightarrow$ Proposition (B)
for $4 \le g \le 9$.}

\medskip

{\bf (2.0)} \ 
Assume that the smooth prime
Fano threefold $X = X_{2g-2}$ ($4 \le g \le 9$)
does not contain a 1-dimensional family of conics of rank 2.

\medskip

{\bf (2.1) Lemma.} \
{\sl
Under the assumption {\bf (2.0)}:

{\bf (i).} \ The Hilbert scheme ${\cal H}_X$ of lines on
$X$ has a unique irreducible component ${\cal H}_o$;

{\bf (ii).} \ ${\cal H}_o = {\cal H}_X$ is exotic,
and the lines $l \in {\cal H}_o$ sweep out
a tangent scroll $R_o \in |{\cal O}_X(d)|, d \ge 2$.
}

\medskip

{\bf Proof of (2.1).} \

{\bf (i).} \
Let ${\cal H}_o$ and ${\cal H}_{\infty}$
be two different irreducible components of ${\cal H}_X$,
and let $R_o$ and $R_{\infty}$ be the surfaces swept
out by the lines $l \in {\cal H}_o$ and $l \in {\cal H}_{\infty}$. 
Since ${\bf Pic}(X) = {\bf Z}.H$, where $H$ is the hyperplane
section, any effective divisor on $X$ must be ample.
In particular the general line $l \in {\cal H}_o$ intersects
the surface $R_{\infty} \subset X$.
If moreover $l \subset R_{\infty}$
for the general $l \in {\cal H}_o$, then $R_o = R_{\infty}$
and this surface contains two 1-dimensional families of lines.
Therefore $R_o = R_{\infty}$ is a quadric surface on $X$,
which contradicts ${\bf Pic}(X) = {\bf Z}.H$.
Therefore the general $l \in {\cal H}_o$ 
intersects $R_{\infty}$ and does not lie in $R_{\infty}$;   
and since $R_{\infty}$ is swept out by lines then
there exists a line $l' \subset R_{\infty}$ which intersects $l$. 
Since $l \in {\cal H}_o$ is general this produces a
1-dimensional family of intersecting lines $l+l'$
(= conics of rank 2 on $X$) 
-- contradiction.

\smallskip

{\bf (ii).} \ 
If ${\cal H}_X = {\cal H}_o$ is non-exotic then,
by Lemma (3.7) of [I1],
the general element of ${\cal H}_o$ will represent a line
$l \subset X$ which will intersect at least one other line
$m \subset X$, i.e. $X$ will contain a 1-dimensional
family of intersecting lines (= conics of rank 2).
Therefore ${\cal H}_o$ is exotic.

Since ${\cal H}_o$ is exotic and $g \ge 4$ then
$R_o$ is the tangent scroll to a curve
$C \subset X$ (see {\bf (1.2)(3)}),
and since ${\bf Pic}(X) = {\bf Z}.H$ then
$R_o \in |{\cal O}_X(d)|$ for some ineteger $d \ge 1$.
If $R_o$ is a hyperplane section of $X$ (i.e. $d = 1$)
then, by Lemma 6 in [Pr],  $C = C_g$ must be a rational
normal curve of degree $g$.
However the last is impossible since then Lemma {\bf (A)}
will imply that $X$ is singular. Therefore $d \ge 2$.
{\bf q.e.d.} 

\smallskip

We shall show that nevertheless
$X$ contains a 1-dimensional family of conics
$l+m$ of rank $2$ where $l,m \in {\cal H}_o$.

\medskip 

{\bf Remark.}
Let $C_s$ be the (possibly empty) set of singular points
of $C$. For any $x \in C - C_s$ denote by $l_x$ the tangent line
to $C$ at $x$. For a point $x = x(0) \in C_s$
define a tangent line to $C$ at $x$ to be any limit
$lim_{x(t) \rightarrow x(0)} \ l_{x(t)}$
of tangent lines $l_{x(t)}$ to points $x(t) \in C - C_s$
(see Ch.2 \S 4 in [GH]).
Clearly, $C$ can have only a finite number of tangent lines to
$x(0) \in C_s$ (see also Ch.2 \S 1.5 in [Sh]).

\medskip 

{\bf (2.2)} \ 
By the initial assumption {\bf (2.0)},
$X$ does not contain a 1-dimensional
family of conics of rank 2.
In particular, $X$ does not contain a 1-dimensional
family of pairs of intersecting tangent lines to $C$.

Since the surface $R_o \in |{\cal O}_X(d)|$,
$d \ge 2$, and $Span \ X = {\bf P}^{g+1}$
then $Span \ R_o = {\bf P}^{g+1}$, $g \ge 4$
(see also Lemma 6 in [Pr]).
Since $R_o$ is swept out by the tangent lines
to $C$ then $Span \ C = Span \ R_o = {\bf P}^{g+1}$.
In particular $C$ does not lie on a plane.
Since $R_o$ is the tangent scroll to the non-plane
curve $C$ then $S_C$ is singular along the curve $C$.

Let $L \not= C$ be (if exists) an irreducible
curve on $R_o$ such that $R_o$ is singular along $L$.
If $L$ is not a tangent line to $C$,
then the general point of $L$ will be an intersection
point of two or more tangent lines to $C$
(see \S 4 in [P2]). The last is
impossible since, by assumption, on $X$ can lie at
most a finite number of pairs of intersecting lines.

Therefore any irreducible curve $L \not= C$
such that $L \subset Sing \ R_o$
must be a tangent line to $C$.
In addition,
the tangent scroll $R_o$ to $C$ still can be singular
along a tangent line $L$ to $C$ --  for example if
$L$ is a common tangent line to two or more
branches of $C$ at $x$, or if $C$ has a branch
with a cusp at $x$, or if $x \in C - C_s$
but $x$ is an inflexion point of $C$ and then $R_o$
has a cusp along $l_x$, etc. (see \S 2, \S 4 in [P2]).

\smallskip

Let $\Delta \subset R_o$ be the union of all the
irreducible curves $L$ on $R_o$ such that
$L \not= C$ and $R_o$ is singular along $L$.
By the above argument, either $\Delta = \emptyset$ or
$\Delta$ is a union of a finite number of tangent
lines to $C$.

Let ${\nu}:R_n \rightarrow R_o$ be the normalization of $R_o$.
Fix a desingularization ${\tau}:\tilde{R}_n \rightarrow R_n$,
and let ${\sigma} = {\tau}\circ{\nu}:\tilde{R}_n \rightarrow R_o$.
Let $E_1,...,E_k$ ($k \ge 0$) be all the irreducible
{\it contractable} curves on $\tilde{R}_n$,
i.e. all the irreducible curves $E_i \subset \tilde{R}_n$
such that ${\sigma}(E_i) \in R_o$ is a point.

Denote by $\sim$ the linear equivalence of divisors
on the smooth surface $\tilde{R}_n$, and let
$E$ be a divisor on $\tilde{R}_n$.
Call $E$ a {\it zero} divisor on $\tilde{R}_n$ if $E \sim 0$;
call the non-zero divisor $E$ {\it contractable}
if $E \sim a_1E_1+...+a_kE_k$
for some $a_1,...,a_k \in {\bf Z}$.
Let $C' \subset \tilde{R}_n$ be the proper
$\sigma$-preimage of $C$ on $\tilde{R}_n$.
Since $R_o$ is the tangent scroll to
the irreducible curve $C$ then
the curve $C'$ is irreducible
and ${\sigma}|_{C'}:C' \rightarrow C$
is an isomorphism over a dense open subset of $C$
(see also Lemma {\bf (2.3)} below).

\smallskip

Let $C'_1,...,C'_r$ be all the irreducible curves
on $\tilde{R}_n$ such that ${\sigma}(C'_i)$
is an irreducible component of
$\Delta$.
%
%
Therefore

\medskip

\centerline{$K_{\tilde{R}_n} \sim {\sigma}^*K_{R_o} - mC'
- {\Sigma}_{i=1,...,r} \ p_iC'_i + E$} 

\medskip
\noindent
for some positive integers $m, p_1,...,p_r$,
and a contractable (or zero) divisor $E$ on $\tilde{R}_n$.

\medskip 

{\bf (2.3) Lemma.}
{\sl
Let $X$ fulfills {\bf (2.0)}. 
Then the tangent scroll $R_o \subset X$ to $C$
has a cusp of type $v^2 = u^3 +...$
along $C$, at a neighbourhood
of the general point $x \in C \subset R_o$;
in particular $m = mult_C \ R_o = 2$ 
} \ 
(see also \S 5 in [P2] and \S 4 in [P1]). 

\medskip

{\bf Proof of (2.3).} \

{\bf (1).} \
We shall see first that $R_o$ is irreducible
at any neighbourhood of the general point $x \in C$,
i.e. $R_o$ has one local branch at $x$.

Assume the contrary, and let $x \in C$ be general.
Let $U \subset X$ be a complex-analytical neighbourhood
of $x$ such that $R_U = R_o \cap U$ is reducible.
Since $R_o$ is swept out by the tangent lines to $C$,
the last imlies that for the general point $y \in C_U = C \cap U$
(hence for the general $y \in C$)
there exists (possibly non-unique) $z \in C$, $z \not= y$ such that
$y$ lies on a tangent line to $C$ at $z$.
Since the set $C_s = \{ x_1,...,x_s \}$ of singular points of $C$
is finite (or empty), and any such $x_i$ has at most a
finite number of tangent lines, then the general $y \in C$ doesn't lie
on a tangent line to $z \in C_s$. Therefore the general $y \in C$
lies on {\it the} tangent line $l_z$ to $C$ at some
(possibly non-unique) $z \in C - C_s$.

If moreover $l_y \not= l_z$ (where $l_y$ is the tangent line
to $C$ at $y$) then all such $l_y + l_z$ will produce a 1-dimensional
family of conics of rank 2 on $X$, which contradicts the
initial assumption {\bf (2.0)} about $X$.

If $l_y = l_z$ then this will imply that the tangent
line $l_y$ to $C$ at the general $y \in C$ is tangent to $C$
at two or more points. But then
the projection $\overline{C}$ of $C \subset {\bf P}^{g+1}$
from the general subspace
${\bf P}^{g-2} \subset {\bf P}^{g+1} = Span \ X$
will be a plane curve with a 1-dimensional family of
lines tangent to $\overline{C}$ at
two or more points, which is impossible.

\smallskip

{\bf (2).} \
It rests to see that the unique local branch of $R_o$
at the general $x \in C$ has a cusp of type $v^2 = u^3 +...$
along $C$ at a neighbourhood of $x$.

Since $R_o \in |{\cal O}_X(d)|$ and $d \ge 2$,
then
$Span \ C = Span \ R_o = {\bf P}^{g+1}$, $g \ge 4$
(see above).

Let $x$ be a general point of $C$.
In order to prove that the tangent scroll
$R_o \subset X$ to $C$ has
a cusp of type $v^2 = u^3 +...$ at a neighbourhood of $x$
it is enough to see that the
projection of $R_o$ from a general
${\bf P}^{g-3} \subset {\bf P}^{g+1}$
has a cusp at $x$.
This reduces the check to the case when $R_o$ is the tangent scroll
to a curve $C \subset {\bf P}^3$.

Since $x$ is a {\it general} point of $C \subset {\bf P}^3$
then, after a possible linear change of coordinates in ${\bf P}^3$,
the curve $C$ has (at $x = (1:0:0:0)$) a local parameterization,
or a {\it normal form} (see \S 2 in [P2], or Ch.2 \S 4 in [GH]):

\medskip

\centerline{
$C_U: \ (x_o(z):...:x_n(z)) =
(1:z + o(z^2): z^2 + o(z^3): z^3 + o(z^4))$, $|z| < 1$,}

\medskip
\noindent
where $o(z^k) = {\Sigma}_{j \ge k} \ a_jz^j$.
Since the coefficient at $z^k$ in $x_k(z) = z^k + o(z^{k+1})$
is $1 \not= 0$ ($k = 2,3$) then,
after a possible linear change of $(x_1,x_2,x_3)$,
the local parameterization of $C$ at
$x = (1:0:0:0)$ can be written as

\medskip

\centerline{
$C_U: \ (x_o(z):x_1(z):x_2(z):x_3(z)) =
(1:z + o(z^4): z^2 + o(z^4): z^3 + o(z^4))$, $|z| < 1$,}

\medskip
\noindent
i.e. $C_U$ approximates, upto $o(z^4)$, the twisted cubic
$C_3 = \{ (1:z:z^2:z^3)$ \}.
Therefore, at a neighbourhood of $x = (1:0:0:0)$,
the unique local branch (see {\bf (1)}) of
the tangent scroll $R_o$ to $C$ is parameterized by

\medskip

\centerline{
$R_U: \ (x_o(z,t):x_1(z,t):x_2(z,t):x_3(z,t))$}

\centerline{
= $(1:z+t+o(z^4)+o(z^3)t:
      z^2+2zt+o(z^4)+o(z^3)t:
      z^3+3z^2t+o(z^4)+o(z^3)t)$.}

\medskip

In affine coordinates $(x_1,x_2,x_3)$ the tangent line to
$C$ at $x = (0,0,0)$ is spanned by the vector $n_x$ = $(1,0,0)$,
and the normal space ${\bf C}^2_o$ $\subset {\bf C}^3(x_1,x_2,x_3)$
to $n_x$ at $x$ is defined by $x_1 = 0$.
In order to prove that $R_U$ has a cusp along $C$
at a neighbourhood of $x$ we shall see that the curve
$D_U$ = $R_U \cap (x_1 = 0)$ $\subset {\bf C}^2_o$
has a cusp at $x$.

\smallskip

On $D_U$ = $R_U \cap (x_1 = 0)$, one has: \
$0 = x_1$ = $z + t + o(z^4) + o(z^3)t$,
i.e. $t$ =  $-z + o(z^4)$.
Let $u= -x_2$, $v = -x_3/2$.
Therefore, on $D_U \subset {\bf C}^2_o$

\medskip

\centerline{
$u = -(z^2 + 2zt + o(z^4) + o(z^3)t) = z^2 + o(z^4)$}

\centerline{
$v = -1/2(z^3+3z^2t+o(z^4)+o(z^3)t) = z^3 + o(z^4))$,}
 
\medskip
\noindent
i.e.
$u^2 = z^4+o(z^6)$,
$uv = z^5+o(z^6)$,
$v^2 = z^6+o(z^7)$,
$u^3 = z^6+o(z^8)$,
$u^2v = z^7+o(z^8)$, ...
%
%
%

Let $C|_U = (f_U(u,v) = 0)$
be the local equation of $D_U \subset {\bf C}^2_o(u,v)$
at $x = (0,0)$.
Therefore, upto a constant non-zero factor,
$f_U(u,v) = v^2 - u^3 + c_{2,1}u^2v + c_{1,2}uv^2 + ...$,
i.e. $C|_U$ has a double cusp-singularity
of type $v^2 = u^3 + ...$ at $x = (0,0)$
(see \S 5 in [P1], \S 4 in [P2],
Ch.5 Examples 3.9.5, 3.9.1 and Ch.1 Exercise 5.14 in [H]).

Therefore $R = R_o$ has a pinch of type $v^2 = u^3 + ... \ $
along $C$ at a neighbourhood of the general point $x \in C$,
which proves Lemma {\bf (2.3)}.

\medskip 

{\bf (2.4)}\ 
By the definition of $C'_1,...,C'_r$
any irreducible component
of $\Delta$ can be represented
(possibly non-uniquely) as the image
${\sigma}(C'_i)$ of some $C'_i, i = 1,...,r$.
Since ${\sigma}|_{C'}:C' \rightarrow C$ is an
isomorphism over an open dense subset of $C$
then the general point $x \in C$ has a unique
${\sigma}$-preimage $x'$ on $C'$,
and the proper preimage $l'_x \subset \tilde{R}_n$
of the tangent line $l_x$ to $C$ at $x$
intersects $C'$ transversally at $x'$.
Since, by assumption, on $X$ doesn't lie a 1-dimensional
family of pairs of intersecting lines
then the tangent line $l_x$ to $C$
at the general point $x \in C$
does not intersect any other tangent line to $C$.

\smallskip

Therefore the non-singular surface $\tilde{R}_n$
has a structure of a possibly non-minimal ruled surface
with a general fiber $L'$ := the proper $\sigma$-preimage
of the general tangent line $l_x$ to $C$.  
In particular $K_{\tilde{R}_n}.L' = -2$, 
and since the curve $C'$ is a section of $\tilde{R}_n$ 
then $C'.L' = 1$. 

By the definition of $C'_i$
the curves ${\sigma}(C'_i) \subset R_o$
are irreducible components of $\Delta$; 
and since by {\bf (2.2)} the components of $\Delta$ 
can be only tangent lines to $C$ 
then ${\sigma}(C'_i)$ is a tangent line to $C$. 
Therefore any component of ${\sigma}^{-1}({\sigma}(C'_i))$, 
in particular $C_i'$, will not intersect the general fiber
$L'$ of $\tilde{R}_n$, i.e. $C'_i.L' = 0$.

Moreover a contractable curve $E_j$
can't intersect the general fiber of $\tilde{R}_n$
since otherwise the point ${\sigma}(E_j) \in R_o$
will be a common point of a $1$-dimensional family
of tangent lines to $C$. 
The last is impossible since $g \ge 4$ and the smooth
$X = X_{2g-2}$ can't contain cones -- see {\bf (1.2)}.
Therefore $E_j.L' = 0$ for any $j = 1,...,k$;
and since $E$ is a sum of such $E_j$ then $E.L' = 0$.

Since $K_X \sim -H$ and $R_o \sim dH$ on $X$
then, by adjunction, $K_{R_o} \sim (d-1)H|_{R_o}$.
Since the hyperplane section $H$ intersects
the general tangent line $l$ to $C$ at one point
then ${\sigma}^*(H|_{R_o})$ is also a section of $\tilde{R}_n$,
i.e. ${\sigma}^*(H|_{R_o}).L' = 1$. 
Therefore 

\medskip 

\centerline{
$-2 = K_{\tilde{R}_n}.L'$
= 
$({\sigma}^*K_{R_o} - mC' - {\Sigma}_{i=1,...,r} \ p_iC'_i + E).L'$
}

\centerline{
=
$(d-1){\sigma}^*(H|_{R_o}).L' - mC'.L'
- {\Sigma}_{i=1,...,r} \ p_iC'_i.L' + E.L'$
=
$(d-1) - m$, i.e. $d = m-1$.
}

\medskip 

Since $X = X_{2g-2}$ is smooth and $g \ge 4$
then, by Lemma {\bf (A)},
$d >1$. Therefore $m = d+1 > 2$,
which is impossible since $m = 2$
by Lemma {\bf (2.3)}.

\medskip

This contradicts the initial assumption {\bf (2.0)}
that the smooth $X = X_{2g-2}$, $4 \le g \le 9$
does not contain a 1-dimensional family of conics of rank 2.
{\bf q.e.d.}

\bigskip

\bigskip

\centerline{\bf \S \ 3. \ Proof of Lemma (A)}

\medskip

{\bf (3.1)} \
By [M1], [M2] any {\it smooth} prime Fano threefold
$X_{2g-2} \subset {\bf P}^{g+1}$, $3 \le g \le 10$
is a complete intersection of hypersurfaces
$F_1,F_2,..., F_N$ of degrees $d_1,d_2,...,d_N$
in a homogeneous (for $g = 6$ -- an almost homogeneous)
space ${\Sigma}(g)$, and:

if $g = 3$ then ${\Sigma}(3) = {\bf P}^4$, $N = 1$,
$d_1 = 4$;

if $g = 4$ then ${\Sigma}(4) = {\bf P}^5$, $N = 2$,
$d_1 = 2$, $d_2 = 3$;

if $g = 5$ then ${\Sigma}(5) = {\bf P}^6$, $N = 3$,
$d_1 = d_2 = d_3 = 2$;

if $g = 6$ then ${\Sigma}(6) = K.G(2,5) \subset {\bf P}^{10}$
is a cone over the grassmannian $G(2,5) \subset {\bf P}^9$, $N = 4$,
$d_1 = d_2 = d_3 = 1, d_4 = 2$;

if $7 \le g \le 10$ then
$X_{2g-2} = {\Sigma}(g) \cap {\bf P}^{g+1}$,
where
${\Sigma}(7) \subset  {\bf P}^{15}$
is the spinor $10$-fold,
${\Sigma}(8) = G(2,6) \subset {\bf P}^{14}$,
${\Sigma}(9) \subset {\bf P}^{13}$ is the sympletic
grassmann $6$-fold,
and ${\Sigma}(10) \subset {\bf P}^{13}$ is the
$G_2$-fivefold.

\medskip 

{\bf (3.2)} \
To prove Lemma {\bf (A)}, it is enough to see that
if $X = X_{2g-2} \subset {\Sigma}(g)$
is a 3-fold complete intersection as in ${\bf (3.1)}$
(assuming implicitly that such $X$ may have singularities)
then $X_{2g-2}$ can't be smooth.
We shall prove this separately for any value of $g$,
$3 \le g \le 9$.

For $g = 3,4,5,6,8$ we use that ${\Sigma}(g)$ is either
a projective space or a (cone over) grassmannian, which makes
it possible to compute directly that the general such
$X_{2g-2} \supset S_{2g-2}$ must have $12-g$ singular points
on the curve $C_g$.

For $g = 7,9$ we assume that $X = X_{2g-2} \subset S_{2g-2}$
is smooth, and then project $X$ from a tangent line to $C_g$
to derive a contradiction on the base of the already known
Lemma {\bf (A)} for $g = 5$.

\medskip 

{\bf (3.3) \ Notation.} \
Let $n \ge 1, m \ge 0$ be integers, let
${\bf P}^{n+m}(z:w) = {\bf P}^{n+m}(z_0:...:z_n:w_{n+1}:...:w_{n+m})$
be the complex projective $(n+m)$-space, and let
$F(z:w) = F(z_0:...:z_n:w_{n+1}:...:w_{n+m})$ be a homogeneous form.
Denote by
${\nabla}_z F =
({\partial F}/{\partial z_0},...,{\partial F}/{\partial z_n })$
the gradient vector of $F$ with respect to $(z) = (z_0:...:z_n)$.

Let $F_1(z:w),...,F_k(z:w)$ ($k \ge 1$) be homogeneous forms.
Denote by:

$(F_1,...,F_k)
\subset {\bf C}[z:w] = {\bf C}[z_0:...:z_n:w_{n+1}:...:w_{n+m}]$
--
the homogeneous ideal generated by $F_1,...,F_k$;

$V(F_1,...,F_k)$
--
the projective variety defined by $F_1,...,F_k$;

$J_z|_{(a:b)}$
=
$J_z(F_1,...,F_k)|_{(a:b)}$
=
$[{\nabla}_z F_1; ...;{\nabla}_z F_k]|_{(a:b)}$
--
the Jacobian matrix $J_z$ of partial derivatives of $F_1,...,F_k$
with respect to $(z)$ = $(z_0,...,z_n)$,
computed at the point $(a:b) \in {\bf P}^{n+m}(z:w)$,
(where ${\nabla}_z F_i$ are regarded as {\it rows} of $J_z$).

Let e.g. $m = 0$, let
$X = V(F_1,...,F_k) \subset {\bf P}^n(z_0:...:z_n)$,
and let $dim \ X = d$.
Then $dim \ T_x \ X \ge d$ for any $x \in X$,
where $T_x \ X$ is the tangent space to $X$ at $x$;
and the point $x \in X$ is singular if $dim \ T_x \ X > d$
(see e.g. Ch. 2, \S 1.4 in [Sh]).
Equivalently $x \in X$ is singular if $rank \ J_z|_x < n-d$.
The subset
$Sing \ X = \{x \in X: rank \ J_z|_x < n-d \} \subset X$
of all the singular points of $X$ is a proper closed subset
of the projective algebraic variety $X$, defined on $X$
by vanishing of all the $(n-d){\times}(n-d)$ minors of $J_z$.

\medskip 

{\bf (3.4) \ Proof of Lemma (A) for g = 3.} \
The tangent scroll to the twisted cubic
$C_3: (x_0:x_1:x_2:x_3) = \vec{s} := (s_0^3:s_0^2s_1:s_0s_1^2:s_1^3)$
is the quartic surface $S_4 = V(f) \subset {\bf P}^3(x)$
where $f(x)$ =
$3x_1^2x_2^2 + 6x_0x_1x_2x_3 - 4x_1^3x_3 -x_0^2x_3^2 -4x_0^2x_2^3$.
The surface $S_4$ is singular along $C_3$
since the gradient vector
${\nabla}_x f|_{\vec{s}} = 0$ for any $\vec{s} \in C_3$.

Let the quartic threefold $X_4 \subset {\bf P}^4(x:u)$
be such that $S_4 = X_4 \cap {\bf P}^3(x)$, and let
$X_4 = V(F) \subset {\bf P}^4(x:u)$ where
$F(x:u) = {\Sigma}_{0 \le k \le 4}\ f_i(x)u^{4-i}$.
Therefore $f_4 \in (f)$, i.e. $f_4 = cf$ for some
constant $c \in {\bf C}$.

Let $x \in X_4$. Then $x \in Sing \ X_4$ {\it iff}
${\nabla}_{x,u} F_x = 0$.
Let $s = (s_0:s_1) \in {\bf P}^1$.
Then $(\vec{s}:0) \in Sing \ X_4$ {\it iff}
$0 = {\nabla}_{x,u} F|_{(\vec{s}:0)}$
=
$({\nabla}_x F,{\partial}F/\partial{u})|_{(\vec{s}:0)}$
=
$({\nabla}_x f|_{\vec{s}}, f_3(\vec{s}))$
=
$(0,f_3(\vec{s}))$.

Therefore
\underline{either}
$f_3(\vec{s}) \equiv 0$
(i.e. $f_3(\vec{s}) = 0$ for any $s = (s_0:s_1)$),
and then $X_4$ is singular along $C_4$,
\underline{or}
$f_s(\vec{s}) \not\equiv 0$,
and then
$(\vec{s}:0) \in Sing \ X_4$
{\it iff}
$s = (s_0:s_1)$ is a zero of the (non-vanishing)
homogeneous form
$F_9(s) = f_3(\vec{s}) = f_3(s_0^3:s_0^2s_1:s_0s_1^2:s_1^3)$
of degree $9$.

\smallskip

In addition, for the general $f_3(x)$
all the zeros of $F_9(s) = f_3(\vec{s})$
are simple, i.e. different from each other.
Therefore the general $X_4 \supset S_4$
has $9 = 12-g(X_4)$ singular points on $C_3$.
In coordinates as above, these singular points
of $X_4$ are the images of the
$9$ zeros of $F_9(s)$ under the Veronese map
$v_3:{\bf P}^1 \rightarrow C_3 \subset X_4$,
$v_3:s = (s_0:s_1) \mapsto (\vec{s}:0)$.

\medskip 

{\bf (3.5) \ Proof of Lemma (A) for g = 4.} \
The tangent scroll to the rational normal quartic
$C_4 : (x_0:x_1:...:x_4) = \vec{s} := (s_0^4:s_0^3s_1:...:s_1^4)$
is a complete intersection
$S_6 = V(q,f) \subset {\bf P}^4(x)$
where $q(x)  =  3x_2^2 - 4x_1x_3 + x_0x_4$
and $f(x)  =  x_2^3 -2x_0x_3^2 -2x_1^2x_4 +3x_0x_2x_4$.
The surface $S_6$ is singular along $C_4$ since the
gradients of $q$ and $f$ are linearly dependent along $C_4$;
more precisely
${\nabla}_x f|_{\vec{s}} = s_0^2s_1^2{\nabla}_x q|_{\vec{s}}$
for any $\vec{s} \in C_4$.

Let $X_6 = V(Q,F) \subset {\bf P}^5(x:u)$ be a complete
intersection of the quadric
$Q(x:z) = {\Sigma}_{0 \le k \le 2} \ q_k(x)u^{2-k} = 0$
and the cubic
$F(x:z) = {\Sigma}_{0 \le l \le 3} \ q_l(x)u^{2-l} = 0$,
and let $S_6 = X_6 \cap {\bf P}^4(x)$.
In particular,
$(q_2,f_3) \subset (q,f)$ as homogeneous ideals in
${\bf C}[x] = {\bf C}[x_0:...:x_4]$.

For the fixed $X_6 = V(Q,F)$ the generators $Q,F$
of the homogeneous ideal $(Q,F)$ can be replaced
by $c'Q$ and $c''F + L(x:u)Q$, for any pair of nonzero
constants $c'$ and $c''$, and for any linear form $L(x:u)$.
Now, $(q_2,f_3) \subset (q,f)$
yields that one can choose $Q$ and $F$ such that
$q_2 = {\varepsilon}'q$ and $f_3 = {\varepsilon}''f$,
where ${\varepsilon}', {\varepsilon}''$ are either $0$ or $1$.

Consider the general case ${\varepsilon}' = {\varepsilon}'' = 1$;
the study in the degenerate case ${\varepsilon}'.{\varepsilon}'' = 0$
is similar. The subscheme $Sing \ X_6  = Sing \ V(Q,F)$ is defined
by
$rank \ [{\nabla}_{x,u} Q ; {\nabla}_{x,u} F] \le 1$.
By the choice of $Q$ and $F$,
${\nabla}_{x,u} Q|_{(\vec{s}:0)} =
({\nabla}_x q_2|_{\vec{s}}, q_1(\vec{s}))$
and ${\nabla}_{x,u} F|_{(\vec{s}:0)} =
({\nabla}_x f_3|_{\vec{s}}, f_2(\vec{s}))$,
where $1.q_2 = q$ and $1.f_3 = f$.
Just as in case $g = 3$, the last together with
the identity
${\nabla}_x f|_{\vec{s}} - s_0^2s_1^2{\nabla}_x q|_{\vec{s}} \equiv 0$
imply that
$(\vec{s}:0) \in Sing \ X_6$
{\it iff}
$F_8(s) := f_2(\vec{s}) - s_0^2s_1^2.q_1(\vec{s}) = 0$,
where $\vec{s} = (s_0^4:s_0^3s_1:...:s_1^4)$.

The Veronese map $v_4:{\bf P}^1 \rightarrow C_4 \subset X_6$,
$v_4:s=(s_0:s_1) \mapsto (\vec{s}:0)$
states an isomorphism between ${\bf P}^1$ and $C_4$.
Therefore \underline{either} $F_8(s) \equiv 0$,
and then $X_6$ is singular along $C_4$,
\underline{or} $F_8(s) \not\equiv 0$,
and then the singular points of $X_6$ on $C_4$
are the $v_4$-images of the zeros of
the homogeneous form $F_8(s)$ of degree $8$.
As in case $g = 3$, for the general
$f_2(x), q_1(x)$ the form $F_8(s)$
has only simple zeros.
Therefore the general
$X_6 \supset S_6$ has $8 = 12-g(X_6)$
singular points on $C_4$.

\medskip 

{\bf (3.6) \ Proof of Lemma (A) for g = 5.} \
The tangent scroll to the rational normal quintic
$C_5 : x_i = s_0^{5-i}s_1^i, 0 \le i \le 5$ is
a complete intersection
$S_8 = V(q',q'',q''') \subset {\bf P}^5(x)$
where
$q'(x)  =  4x_1x_3-3x_2^2-x_0x_4$,
$q''(x)  =  3x_1x_4-2x_2x_3-x_0x_5$
and
$q'''(x)  =  x_1x_5-4x_2x_4+3x_3^2$.
The surface $S_8$ is singular along $C_4$ since the
gradients of $q'$, $q''$ and $q'''$
are linearly dependent along $C_5$. More precisely

\medskip

$\bf ({\ast})$ \ \ \ \ \ \ \ \
$s_1^2{\nabla}_x q'|_{\vec{s}} - s_0s_1{\nabla}_x q''|_{\vec{s}} -
s_0^2{\nabla}_x  q'''|_{\vec{s}} = 0$ for any $\vec{s} \in C_5$.

\medskip

Let $X_8 = V(Q',Q'',Q''') \subset {\bf P}^6(x:u)$ be a complete
intersection of the quadrics
$Q^i(x:z) = {\Sigma}_{0 \le k \le 2} \ q^i_k(x)u^{2-k}$
($i = ','','''$), and such that
$S_8 = X_8 \cap {\bf P}^5(x)$.
In particular $(q'_1,q''_1,q'''_1) \subset (q',q'',q''')$
as homogeneous ideals in ${\bf C}[x] = {\bf C}[x_0:...:x_4]$.
Therefore $q'_1$, $q''_1$ and $q'''_1$ are linear combinations
of $q'$, $q''$ and $q'''$. Since the ideal
$(Q',Q'',Q'')$ of $X_8$ is generated also by
any $GL(3)$-transform of the triple $(Q',Q'',Q''')$,
we may assume that
$q'_2 = {\varepsilon}'q'$, $q''_2 = {\varepsilon}''q''$
and
$q'''_2 = {\varepsilon}'''q'''$, where ${\varepsilon}'$, ${\varepsilon}''$
and ${\varepsilon}'''$ are $0$ or $1$.
The subscheme $(Sing \ X_8)|_{C_5} \subset {\bf P}^1$ is defined by

\medskip

$\bf ({\ast \ast})$ \ \ \ \ \ \ \ \
$2 \ge rank \ [{\nabla}_{x,u} Q';\
{\nabla}_{x,u} Q'';\ {\nabla}_{x,u} Q''']|_{(\vec{s}:0)}$

\centerline{= $rank \
[({\nabla}_x {\varepsilon}'q'|_{\vec{s}}, q'_1(\vec{s}));\
({\nabla}_x {\varepsilon}''q''|_{\vec{s}}, q''_1(\vec{s}));\
({\nabla}_x {\varepsilon}'''q'''|_{\vec{s}}, q'''_1(\vec{s})]$.}

\medskip

Let
$F_7(s) =
{\varepsilon}'s_1^2q'_1(\vec{s}) -
{\varepsilon}''s_0s_1q''_1(\vec{s}) -
{\varepsilon}'''s_0^2q'''_1(\vec{s})$.

The Veronese map $v_5:{\bf P}^1 \rightarrow C_5 \subset X_8$,
$v_5:s = (s_0:s_1) \mapsto (\vec{s}:0)$
states an isomorphism between ${\bf P}^1$ and $C_5$.
Just as in {\bf (3.5}), \ $\bf ({\ast})$ and $\bf ({\ast \ast})$
imply that
\underline{either} $F_7(s) \equiv 0$,
and then $X_8$ is singular along $C_5$,
\underline{or} $F_7(s) \not\equiv 0$,
and then the singular points of $X_6$ on $C_5$
are the $v_5$-images of the zeros of
the homogeneous form $F_7(s)$ of degree $7$.
Moreover, for the general linear forms
$q'_1(x),q''_1(x),q'''_1(x)$
the form $F_7(s)$ has only simple zeros.
Therefore the general
$X_8 \supset S_8$ has $7 = 12-g(X_8)$
singular points on $C_5$.

\bigskip

\bigskip

\centerline{\bf Lemma (A) for g = 6,8.}

\medskip

{\bf (3.7) \ Lemma.}
{\sl
Let $n \ge 3$, and let
$G = G(1:n) = G(2,n+1) \subset {\bf P}^{n(n+1)/2-1}$
be the grassmannian of lines in
${\bf P}^n = {\bf P}({\bf C}^{n+1})$.
Let $C \subset G \subset {\bf P}^{n(n+1)/2-1}$
be a smooth irreducible curve
such that $dim \ Span(C) \ge 3$.
Let the surface $Gr(C) \subset {\bf P}^n$ be the union of lines
$l \subset {\bf P}^n$ such that $l \in C \subset G = G(1:n)$,
and let $S_C \subset Span(C) \subset {\bf P}^{n(n+1)/2-1}$
be the surface swept out by the tangent lines to $C$
(or the tangent scroll to $C$ -- see above).

Then the tangent scroll $S_C$ to $C$ lies on $G = G(1:n)$
{\it iff}

\underline{either}
all the lines $l \in C$ have a common point,
i.e. $Gr(C)$ is a cone,

\underline{or}
there exists an irreducible curve $Z \subset {\bf P}^n$
such that all the lines $l \in C$ are tangent lines to $Z$,
i.e. $Gr(C)$ is the tangent scroll to $Z$.
}

\medskip

{\bf Proof.}
For $n=3$ this result can be found in [AS]. For $n>3$ one
can apply induction, using the fact that projection from a point in
${\bf P}^n$ induces a projection from $G(2,n+1)$ onto $G(2,n)$.

\medskip

{\bf (3.8) \ Lemma.}
{\sl
Let $C_g \subset G(2,g/2+2) = G(1:{\bf P}^{g/2+1})$
($g = 6,8$) be a rational normal curve such that
the tangent scroll $S_{2g-2}$ to $C_g$ is contained
in $G(2,g/2+2)$. Then:

{\bf (i).} \ If $g = 6$ then the lines $l \subset {\bf P}^4$,
$l \in C_6$ sweep out the tangent scroll to a rational
normal curve $C_4 \subset {\bf P}^4$.

{\bf (ii).} \ If $g = 8$, and if there exists a
3-fold linear section $X_{14} = G(2,6) \cap {\bf P}^9$
such that $X_{14} \supset S_{14}$, then the lines
$\{ l \subset {\bf P}^5: l \in C_8 \}$
sweep out the tangent scroll to a rational
normal curve $C_5 \subset {\bf P}^5$.
}

\medskip

{\bf Proof.}
Let $g = 6,8$, let the lines $l \in C_g$ sweep out a cone
$Gr(C_g) \subset {\bf P}^{g/2+1}$ -- see {\bf (3.7)},
and let $x$ be the vertex of $Gr(C_g)$.
Then $C_g$ is contained in the Schubert $g/2$-space
${\bf P}^{g/2}_x = {\sigma}_{g/2,0}(x)$ =
$\{ l \subset {\bf P}^{g/2+1}: x \in l \}$.
Since $C_g$ is projectively normal it must span
a $g$-space. Therefore $g \le g/2+1$ which contradicts
$g = 6,8$. Therefore, by {\bf (3.7)}, the lines $l \in C_g$
must sweep out the tangent scroll to a rational curve
$C_d \subset {\bf P}^{g/2+1}$.

Let $g = 6$, and let $C_d \subset {\bf P}^3 \subset {\bf P}^4$.
Then $C_6 \subset G(1:{\bf P}^3)$ =
${\sigma}_{11}({\bf P}^3) \subset G(1:{\bf P}^4)$.
Therefore $6 = dim(Span \ C_6) \le dim(Span \ G(1:{\bf P}^3)) = 5$
-- contradiction. Therefore $d \ge 4$ since $C_d$ must span
${\bf P}^4$, and now it is easy to see that the rational normal
curve $C_6$ is the curve of tangent lines to $C_d$ {\it iff} $d = 4$.

Let $g = 8$, and let $C_d \subset {\bf P}^4 \subset {\bf P}^5$.
Then $C_8 \subset G(1:{\bf P}^4)$ =
${\sigma}_{11}({\bf P}^4) \subset G(1:{\bf P}^5)$.
Let ${\bf P}^9_o = Span \ G(1:{\bf P}^4)$.
Then ${\bf P}^8 = Span \ C_8 \subset{\bf P}^9_o$.
By condition $C_8 \subset X_{14} = G(1:{\bf P}^5) \cap {\bf P}^9$.
Therefore ${\bf P}^8 \subset {\bf P}^9 \cap {\bf P}^9_o$
and $X_{14} \supset Z_o := G(1:{\bf P}^4) \cap {\bf P}^8$.
Since ${\bf P}^8 \subset {\bf P}^9_o = Span \ G(1:{\bf P}^4)$,
$Z_o$ is at least a hyperplane section of the 6-dimensional
grassmannian $G(1:{\bf P}^4)$. This contradicts $X_{14} \supset Z_o$
and $dim \ X_{14} = 3$.
Therefore $d \ge 5$ since $Span \ C_d = {\bf P}^5$,
and now it is easy to see that the rational normal curve $C_8$
is the curve of tangent lines to $C_d$ {\it iff} $d = 5$.
{\bf q.e.d.}

\bigskip

\centerline{\bf Proof of Lemma (A) for g = 6.}

\medskip

{\bf (3.9)} \
By ${\bf (3.1)}$ any smooth prime $X_{10}$ is
a complete intersection of three hyperplanes and a quadric
in the cone $K.G(2,5) \subset {\bf P}^{10}$.
Let $o$ be the vertex of the cone $K.G(2,5) \subset {\bf P}^{10}$,
and let $X \subset K.G(2,5)$, be as in ${\bf (3.1)}$.
There are two kinds of such threefolds $X_{10}$
(see [I1], [Gu]):

{\bf (i).} \
{\bf g = 6 -- first kind}: \ $o \not\in Span(X)$, and then
the projection $p_o$ from $o$ sends $X$ isomorphically to
$X_{10} = G(2,5) \cap {\bf P}^7 \cap Q$, where
$G(2,5) \subset {\bf P}^9$ by the Pl\"ucker embedding,
${\bf P}^7 \subset {\bf P}^9$ and $Q$ is a quadric.

{\bf (ii).} \
{\bf g = 6 -- second kind}: \ $o \in Span(X)$, and then
${\pi} = {p_o}|_{X}:X = X_{10}' \rightarrow Y_5$ is a double
covering of a threefold $Y_5 = G(2,5) \cap {\bf P}^6$.
In particular, if $X_{10}'$ is smooth then the intersection
$Y_5 = G(2,5) \cap {\bf P}^6$ is smooth.

\bigskip

\centerline{\bf g = 6 (first kind)}

\medskip

{\bf (3.10)} \
Let $X_{10} = G(2,5) \cap {\bf P}^7 \cap Q$
be a (possibly singular) complete intersection
as in {\bf (3.9)(i)},
and assume that $X_{10}$ contains the tangent scroll
$S_{10}$ to a rational normal curve $C_6$ of degree 6.
By Lemma {\bf (3.8)(i)} the points of $C_6$ are the
Pl\"ucker coordinates $x_{ij}(s)$ of the tangent lines
to a rational normal quartic
$C_4 = x_i = s^i, 0 \le i \le 4$, i.e.
\[
(x_{ij}(s)) =
\left(
\begin{array}{ccccc}
 0  & 1   & 2s  & 3s^2 & 4s^3 \\
... & 0   & s^2 & 2s^3 & 3s^4 \\
... & ... & 0   & s^4  & 2s^5 \\
... & ... & ... & 0    & s^6  \\
... & ... & ... & ...  & 0
\end{array}
\right)
\]
Therefore the subspace ${\bf P}^6 = Span(C_6) \subset {\bf P}^9$
is defined by
$H_0 = H_1 = H_2 = 0$ where
$H_0 = x_{03} - 3x_{12}$,
$H_1 = x_{04} - 2x_{13}$
and
$H_2 = x_{14} - 3x_{23}$.

\medskip 

{\bf (3.11)} \
Introduce in ${\bf P}^6$ the coordinates
$(v) = (v_0:...:v_6) =
(x_{01}:x_{02}:x_{12}:x_{13}:x_{23}:x_{24}:x_{34})$,
and let $I_{S_{10}} \subset {\bf C}[v] := {\bf C}[v_0,...,v_6]$
be the homogeneous ideal of the
tangent scroll $S_{10} \subset {\bf P}^6$ to $C_6$.
Let
\
$Pf_k = x_{ab}x_{cd} - x_{ac}x_{bd} + x_{ad}x_{bc}$,
$0 \le a < b < c < d \le 4$: $a,b,c,d \not= k$
$(k \in \{0,1,2,3,4\})$
\
be the $5$ Pl\"ucker quadrics in the coordinates $x_{jk}$.
In coordinates $(v)$ of ${\bf P}^6 = {\bf P}^6(v)$
the restrictions $q_k$ of $Pf_k$ to ${\bf P}^6$ are

\medskip

\centerline{$q_0(v)$ = $v_2v_6 - v_3v_5 + 3v_4^2$, \
$q_1(v)$ = $v_1v_6 - 3v_2v_5 +2v_3v_4$, \
$q_2(v)$ = $v_0v_6 - 9v_2v_4 + 2v_3^2$,}

\centerline{$q_3(v)$ = $v_0v_5 - 3v_1v_4 +2v_2v_3$, \
$q_4(v)$ = $v_0v_4 - v_1v_3 + 3v_2^2$.}

\medskip

In the open subset $U_0 = \{v_0 = 1\} \subset {\bf P}^6(v)$
the curve $C_6$ is parameterized by
$C_6 = \{ (v) = \vec{s} := (1,2s,s^2,2s^3,s^4,2s^5,s^6)\}$.

Now, it is easy to see that the quadric
$q(v) = 5v_2v_4 - 2v_1v_5 +3v_0v_6$ vanishes at
the points of the tangent scroll $S_{10}$ to $C_6$,
and the homogeneous ideal
$I_{S_{10}} = (q_0,...,q_4,q) \subset {\bf C}[v]$.

Let $J_v = J_v(q_0,...,q_4,q)$ be the Jacobian matrix
of $(q_0,...,q_4,q)$ with respect to $v = (v_1,...,v_6)$.
The surface $S_{10} \subset V_5$ is singular along $C_6$
since
$rank \ J_v|_{\vec{s}} < 4 = codim_{{\scriptstyle {\bf P}^6}} \ S_{10}$
for any $\vec{s} \in C_6$. For the special choice $(v)$ of the coordinates
this can be verified directly:

\medskip 

\centerline{$s^{-2}{\nabla}_v q_0|_{\vec{s}} = (0,s^4,-2s^3,6s^2,-2s,1)$,
$s^{-1}{\nabla}_v q_1|_{\vec{s}} = (s^5,-6s^4,2s^3,4s^2,-3s,2)$,}

\centerline{${\nabla}_v q_2|_{\vec{s}} = (0,-9s^4,8s^3,-9s^2,0,1)$,
$s{\nabla}_v q_3|_{\vec{s}} = (-3s^5,4s^4,2s^3,-6s^2,s,0)$,}

\centerline{$s^2{\nabla}_v q_4|_{\vec{s}} = (-2s^3,6s^4,-2s^3,s^2,0,0)$,
and
${\nabla}_v q|_{\vec{s}} = (-4s^5,5s^4,0,5s^2,-4s,3)$.}

\medskip

${\bf ({\ast})}$. \ Therefore \
$rank \ J_v|_{\vec{s}} = 4$ for any $\vec{s} \in C_6$,
and the linear 4-space of linear equations between the gradients
${\nabla}_v q_0|_{\vec{s}},...,{\nabla}_v q_4|_{\vec{s}}$
and ${\nabla}_v q|_{\vec{s}}$
is spanned on the Pfaff-equation
$s^{-2}{\nabla}_v q_0|_{\vec{s}}
+ s^{-1}{\nabla}_v q_1|_{\vec{s}}
+ {\nabla}_v q_2|_{\vec{s}}
+ s{\nabla}_v q_3|_{\vec{s}}
+ s^2{\nabla}_v q_4|_{\vec{s}}= 0$
and the 3-space of equations
${\nabla}_v q|_{\vec{s}}
= a_0.s^{-2}{\nabla}_v q_0|_{\vec{s}}
+ a_1.s^{-1}{\nabla}_v q_1|_{\vec{s}}
+ a_2.{\nabla}_v q_2|_{\vec{s}}
+ a_3.s{\nabla}_v q_3|_{\vec{s}}
+ a_4. s^2{\nabla}_v q_4|_{\vec{s}} = 0$ where:

\medskip

${\bf ({\ast \ast})}$. \ \ \ \
$a_0 + 4a_3 + 3a_4 = 8$,
$a_1 - 3a_3 - 2a_4 = -4$,
$a_2 + 2a_3 + a_4 = 3$.

\medskip 

{\bf (3.12)} \
In the dual space $\check{\bf P}^9$,
let $\check{G} = \check{G}(2,5) \subset \check{\bf P}^9$ be the
grassmannian of hyperplane equations represented
by the skew-symmetric $5 \times 5$ matrices of rank $2$.
It is easy to see that the plane
${\Pi} = <H_0,H_1,H_2> \subset \check{\bf P}^9$
of hyperplane equations of ${\bf P}^6 = Span \ C_6$
does not intersect $\check{G}$.

Let $\Lambda \subset \Pi$ be any line in $\Pi$.
In ${\bf P}^9$, the line $\Lambda$ defines, by duality,
the subspace ${\bf P}^7({\Lambda}) \supset {\bf P}^6 = Span \ C_6$.
It is easy to see that the fourfold
$W({\Lambda}) = G \cap {\bf P}^7({\Lambda})$,
where $G = G(2,5)$, is smooth.
In fact, $W({\Lambda})$ will be smooth
{\it iff} the line $\Lambda$ does not intersect $\check{G}$.
The last is true since ${\Lambda} \subset {\Pi}$ and
${\Pi} \cap \check{G} = \emptyset$. Therefore any
$X_{10} \supset S_{10}$ is a quadratic section of the smooth
4-fold $W = W({\Lambda})$.

\medskip 

{\bf (3.13)} \
Let $X_{10} = G(2,5) \cap {\bf P}^7 \cap Q \supset S_{10}$,
where $Q$ is a quadric.
We shall show that the singularities of $X_{10}$ on $C_6$
are the zeros of a homogeneous form of degree $6$
on $C_6 \cong {\bf P}^1$.
To simplify the notation, we shall show this for one special
choice of the line $\Lambda \subset \Pi$ (see {\bf (3.12)};
the check for any other $\Lambda \subset \Pi$ is similar.

Let $\Lambda = \{ H_2 = 0 \} \subset \Pi$. Then, in coordinates
$(v)$ and $u = x_{14}-3x_{23}$ in ${\bf P}^7({\Lambda})$,
the subspace ${\bf P}^6(x) = (u = 0)$.
Therefore any quadric $Q \subset {\bf P}^7({\Lambda})$,
such that $Q \cap X_{10} = S_{10}$,
can be written in the form $Q = Q(v, u) = cu^2 + L(v)u + q(v)$
where $c \in {\bf C}$ and $L$ is a linear form of $(v)$.

Let $Q_k(v,u)$ ($k = 0,1,...,4$) be the restriction of the
Pfaff quadric $Pf_k$ on ${\bf P}^7(v,u)$, and let
$J_{v,u} =
[{\nabla}_{v,u} Q_0 ; ... ; {\nabla}_{v,u} Q_4 ; {\nabla}_{v,u} Q_0]$
be the Jacobian matrix of $(Q_0,...,Q_4,Q)$.

The singularities of $X_{10}$ on $C_6$ are the points
$(\vec{s}:0) \in {\bf P}^7$ for which
$rank \ J_{v,u}|_{(\vec{s}:0)} < 4$ = $codim_{{\bf P}^7} \ X_{10}$.

Let $l_i(\vec{s}) = {\partial Q_i}/{\partial u}|_{\vec{s}}$,
$i = 0,...,4$.
The rows of $J_{v,u}|_{(\vec{s},0)}$
are
${\nabla}_{u,v} Q_i|_{(\vec{s}:0)} =
({\nabla}_u q_i|_{\vec{s}}, l_i(\vec{s}))$,
$i = 0,...,4$
and
${\nabla}_{u,v} Q|_{(\vec{s}:0)} =
({\nabla}_u q|_{\vec{s}}, L(\vec{s}))$.
For the special choice of the line $\Lambda \subset \Pi$,
the linear forms $l_i = l_i(\vec{s})$ are
$(l_0,l_1,l_2,l_3,l_4) = (s^4, 0, -3s^2, -2s, 0)$.

Therefore, in view of ${\bf ({\ast})}$,
$X_{10}$ will have a singularity at $(\vec{s},0) \in C_6$
if there exist constants $a_0,...,a_4$ satisfying
${\bf ({\ast \ast})}$ and such that
$L(\vec{s})$ = $(a_0 - 3a_2 - 2a_3)s^2$
(here $\vec{s}$ = $(1,2s,s^2,2s^3,s^4,2s^5,s^6)$ -- see above).

\smallskip

By ${\bf ({\ast \ast})}$ \ $a_0-3a_2-2a_3 = -1$.
Therefore, in homogeneous
coordinates $(s_0:s_1)$, $s = s_1/s_0$,
the variety $X_{10}$ will be singular
at the point $(\vec(s):0) \in C_6$
{\it iff}
$F_6(s_0:s_1) :=
L(s_0^6:2s_0^5s_1:s_0^4s_1^2:2s_0^3s_1^3:s_0^2s_1^4:2s_0s_1^5:s_1^6)
+ s_0^4s_1^2 = 0$.

The Veronese map $v_6:{\bf P}^1 \rightarrow C_6 \subset X_{10}$,
$v_6:(s_0:s_1) \mapsto (\vec{s}:0)$
=
$(s_0^6:2s_0^5s_1:s_0^4s_1^2:2s_0^3s_1^3:s_0^2s_1^4:2s_0s_1^5:s_1^6:0)$,
states an isomorphism between ${\bf P}^1$ and $C_6$.
Therefore
\underline{either} $F_6(s_0:s_1) \equiv 0$,
and then $X_{10}$ is singular along $C_6$,
\underline{or} $F_6(s_0:s_1) \not\equiv 0$,
and then the singular points of $X_{10}$ on $C_6$
are the $v_6$-images of the zeros of
the homogeneous form $F_6(s_0:s_1)$ of degree $6$.

\smallskip

The choice of an arbitrary line ${\Lambda} \subset {\Pi}$
will only change the homogeneous sextic forms defined by
$l_i(\vec{s}), i = 0,...,4$. \ {\bf q.e.d.}

\bigskip

\centerline{\bf g = 6 (second kind)}

\medskip

{\bf (3.14)} \
Let $\pi : X = X_{10}' \rightarrow Y_5$, $X \subset {\bf P}^7$,
be a double covering of the Del Pezzo threefold
$Y_5 = G(2,5) \cap {\bf P}^6$
branched along the quadratic section $B \subset Y_5$.
Below we shall identify the branch locus $B \subset Y_5$
and the ramification divisor $R \subset X_{10}'$, $R \cong B$.

\medskip 

{\bf (3.15)} \
Assume that
$X$ contains the tangent scroll $S = S_{10}$
to the rational normal sextic $C = C_6$,
and let
$l \subset S$ be a general tangent line to $C$.
Then
$N_{l/X} \cong {\cal O}(1) \oplus {\cal O}(-2)$
(see {\bf (1.2)}),
and ${\pi}(l) \subset Y_5$ also is a line.
If
$N_{{\pi}(l)/Y_5} \cong {\cal O} \oplus {\cal O}$
then
$d{\pi}: N_{l/X} \rightarrow N_{{\pi}(l)/Y_5}$
has one-dimensional kernel along $l$.
In this case
$l$ is contained in the ramification divisor $R$ of ${\pi}$.
The last is possible only for a finite number of $l's$.
Therefore,
at least for the general tangent line $l \subset S$
to $C$, the line ${\pi}(l) \subset {\pi}(S)$
is a $(-1,1)$-line on $Y_5$,
i.e.
$N_{{\pi}(l)/Y_5} \cong {\cal O}(-1) \oplus {\cal O}(1)$.
Therefore
${\pi}(S) \subset Y_5$
coincides with the surface $S_{-1,1} \subset Y_5$
swept-out by the $(-1,1)$-lines on $Y_5$,
which in turn
is a tangent scroll to a rational normal
sextic (see [FN]). Clearly ${\pi}:S \rightarrow S_{-1,1}$
is an isomorphism and $S_{-1,1}$ is the tangent scroll
to ${\pi}(C_6)$.

Assume first that $S = R$.
Then $S_{-1,1} = {\pi}(S) = B \cong R$.
Therefore $X$ is singular along $C_6$
since the branch locus $B = S_{-1,1}$ of ${\pi}$
is singular along $C_6$.
In order to prove Lemma {\bf (A)}
for $g=6$ (second kind) it rests
to see the general $X_{10}' \supset S_{10}$
such that $S \not= R$ is singular.
This will imply that any
$X_{10}' \subset S_{10}$ must
be singular -- since the property $X_{10}' \subset S_{10}$
to have a singularity is a closed condition.

Since
${\pi}:X_{10}' \rightarrow Y_5$ is a 2-sheeted covering
with a ramification divisor $R$, then
$H^0(X_{10}',O(1)) = {\pi}^*H^0(Y_5,O(1)) + {\bf C}.R$
where $R$ is the hyperplane equation of $R \subset X_{10}'$.
Therefore,
since $S \cong {\pi}(S)$ and $S \not= R$,
then ${\pi}(S)$ is tangent to the branch locus $B$
of $\pi$ along a (possibly singular) canonical
curve $C^6_{10} = B \cap H = {\pi}(S) \cap H$
for some hyperplane $H \subset {\bf P}^6$.
Since ${\pi}(S) = S_{-1,1}$ then the general
$X_{10}' \supset S = S_{10}$, $S \not= R$
comes from a branch locus $B$ totally tangent to
$S_{-1,1}$ along a general hyperplane section
$H \cap S_{-1,1}$.
Since $H$ is general then $H$ intersects the rational
normal sextic ${\pi}(C_6)$ at $6$ points $p_1,...,p_6$
such that $p_i \not= p_j$ for $i \not= j$.
By $p_i \in B$ and the identification $B = R$, we may cosider
$p_i$ as points on $R \subset X_{10}'$. We shall see that

\medskip

{\bf $({\bf \ast})$ Lemma}. \ $p_i \in Sing \ X_{10}'$, $i = 1,...,6$.

\medskip

{\bf Proof of $({\bf \ast})$}. \
Let $U_i \subset Y_5$ be a sufficiently small neighborhood of
the point $p_i$.
Since $Y_5$ is smooth, we can identify $U_i$ with a disk in
${\bf C}^3$, and let $(u,v,w)$ be local coordinates in $U_i$
s.t. $p_i = (0,0,0)$.
Let $f(u,v,w) = 0$ be the local equation of $S_{10}$ in $U_i$.
Since $S_{10}$ is the tangent scroll to the smooth rational
curve $C_{6}$, and $p_i \in C_6$, one can choose the coordinates
$u,v,w$ such that $f = u^3 - v^2$ (since the scroll $S_{10}$
has a double cusp-singularity along $C_{6}$
-- see Lemma {\bf (2.3)}).
Let $q = 0$ and $h = 0$ be the local equations of $B$ and $H$
in $U_i$. Since $B$ and $S_{10}$ are singular along
$C^6_{10}$ = $B \cap H$ = $S_{10} \cap H$, \ $h^2$, $f$ and $q$
are linearly dependent, i.e.
${\alpha}h^2 + {\beta}f + {\gamma}q = 0$ for some constants
${\alpha}, {\beta}, {\gamma}$.
Moreover ${\alpha} \not= 0$ and ${\gamma} \not= 0$
since ${\pi}(S) \cong S$ and ${\pi}(S) \not= B$,
and one may assume that ${\beta} \not= 0$
(otherwise $B = 2H$ and $X_{10}'$ will be singular along the
surface $R_{red} \cong B_{red} = H$).
Since $p_i \in H$, $h(0,0,0) = 0$ and
$h = au + bv + cw + o(2)$, where $o(k)$ denotes a sum of terms
of degree $\ge k$.
Therefore the surface $B$ is singular at $p_i$ since
$B \cap U_i = (q = 0)$, and the series expansion

\medskip

${}$\hspace{1cm}{$q(u,v,w) = -{\beta}/{\gamma}f - {\alpha}/{\gamma}h^2$
= ${\beta}/{\gamma}(v^2-u^3) -{\alpha}/{\gamma}(au+bv+cw + o(2))^2$}
 
\medskip
\noindent
has no linear term. Since $B$ is the branch locus of
${\pi}:X_{10}' \rightarrow Y_5$, the threefold $X_{10}'$
will be also singular at $p_i \in R \cong B$,
$i = 1,...,6$.
In addition, the $6 = 12 - g(X_{10}')$ singular points
$p_1,...,p_6$ of $X_{10}'$ lie on the rational normal
curve $C_6 \subset X_{10}'$. {\bf q.e.d.}

\bigskip

\bigskip

\centerline{\bf Proof of Lemma {\bf (A)} for g = 8.}

\medskip

{\bf (3.16) \ The Da Palatini construction} (see [Pu]). \
Let ${\bf P}^5 = {\bf P}(V)$ where $V = {\bf C}^6$,
and let $\hat{V} = Hom_{\bf C}(V,{\bf C})$ be the dual space of $V$.
The points
$H \in {\wedge}^2 \hat{V} = Hom_{\bf C}({\wedge}^2 V, {\bf C})$
can be regarded as skew-symmetric linear maps
$H:V \rightarrow \hat{V}$,
and the hyperplanes
$(H = 0) \subset {\bf P}^{14} = {\bf P}({\wedge}^2 V)$
can be regarded as points of
$\hat{\bf P}^{14} = {\bf P}({\wedge}^2 \hat{V})$.

Let $Pf  = \{H \in {\wedge}^2 \hat{V}): rank(H) \le 4 \}/{\bf C}^*$.
Then $Pf$ is the Pfaff cubic hypersurface in $\hat{\bf P}^{14}$
defined, in coordinates, by vanishing of the cubic Pfaffian
of the skew-symmetric $(6 \times 6)$-matrix $H$.
Let $U_{10} \subset {\wedge}^2 V$ be a 10-dimensional subspace,
and let

\medskip

\centerline{$\hat{U}_5$ = $U_{10}^{\perp}$
=
$\{ H \in {\wedge}^2 \hat{V}$ =
$Hom_{\bf C}({\wedge}^2 V,{\bf C}): H|_{U_{10}} = 0 \}$.}

\medskip

Let moreover $U_{10} \subset {\wedge}^2 V$ be such that
$rank(H) \ge 4$ for any $H \in U_{10}^{\perp}$,
and let $X_{14}$ = $G(2,6) \cap {\bf P}(U_{10})$
and $B_3$ = $Pf \cap {\bf P}(\hat{U}_5)$.

The construction ``Da Palatini'' of G. Fano
shows that any hyperplane ${\bf P}^4 \subset {\bf P}(V)$
defines a birational isomorphism
${\xi}: X_{14} \rightarrow B_3$
which can be described as follows (see [Pu]):

\smallskip

Identify the point $b \in B_3$
and the (projective equivalence class of) the skew-symmetric
$6 \times 6$ matrix corresponding to $b$.
Since $rank \ b = 4$ for any $b \in B_3$
the projective kernel $n_b$ of $b$ will be a line
in ${\bf P}^5$. Let
$W$ := $\cup_{b \in B_3} ~  \{n_b = {\bf P}(Ker \ b) \} \subset {\bf P}^5$.

The fourfold $W$ can be described by an alternative way.
Identify the point $l \in G(1:5)$ and the line $l \subset {\bf P}^5$,
and let
$W' := \cup \{ l \subset {\bf P}^5 | l \in X_{14} \} \subset {\bf P}^5$.
Then (see [Pu, p. 83]):

\medskip

{\it
{\bf (a).} \ for the general $v \in W'$ there exists a unique
$l \in X_{14}$  such that $v \in l$;

{\bf (b).} \ for the general $w \in W$ there exists a unique
$b \in B$ such that $w \in n_b$;

{\bf (c).} \ $W' = W$.
}

\medskip

Let $H \subset {\bf P}^5$ be a general hyperplane.
Then, by {\bf (a)}, {\bf (b)} and {\bf (c)}, the maps
${\phi}: X_{14} \rightarrow H \cap V$, ${\phi}(l) = H \cap l$
\ and \
${\psi}: B_3 \rightarrow H \cap V$, ${\psi}(b) = H \cap n_b$
are birational isomorphisms.
The composition
${\xi} = {\xi}_H = {\psi}^{-1} \circ {\phi}: X_{14} \rightarrow B_3$ is
a birational isomorphism,
depending on the choice of the hyperplane $H \subset {\bf P}^5$
(see [Pu, p. 85]).

\medskip 

{\bf (3.17) \ The dual cubic fourfold of $S_{14}$.} \
By Lemma {\bf (3.8)(ii)} the curve $C_8$ is the Pl\"ucker image
of the tangent scroll to a rational normal quintic
$C_5: (x_0:...:x_5) = \vec{s} = (1:s:s^2:...:s^5)$
in  ${\bf P}^5(x)$.
The points of the curve $C_8$ are the Pl\"ucker coordinates
$x_{ij}(s)$ of the point $\vec{s} \in C_5$:
\[
(x_{ij}(s)) =
\left(
\begin{array}{cccccc}
 0  &  1  & 2s   & 3s^2 & 4s^3 & 5s^4  \\
... &  0  &  s^2 & 2s^3 & 3s^4 & 4s^5  \\
... & ... &  0   &  s^4 & 2s^5 & 3s^6  \\
... & ... & ...  &  0   &  s^6 & 2s^7  \\
... & ... & ...  & ...  &  0   &  s^8  \\
... & ... & ...  & ...  & ...  &  0
\end{array}
\right)
\]
Therefore $S_{14} = G(2,6) \cap {\bf P}^8$ where
${\bf P}^8 = (H_0 = ... = H_5 = 0) \subset {\bf P}^{14}$ and:

\medskip

\centerline{$H_0 =  x_{03} - 3x_{12}$,
$H_1 =  x_{04} - 2x_{13}$,
$H_2 = 3x_{05} - 5x_{14}$,}

\centerline{$H_3 =  x_{14} - 3x_{23}$,
$H_4 =  x_{15} - 2x_{24}$,
$H_5 =  x_{25} - 3x_{34}$.}

\medskip

Let ${\Pi}^5 := <H_0,...,H_5> \subset \hat{\bf P}^{14}$
be the projective linear span of $\{ H_0,...,H_5 \}$,
and let ${\bf B}_3 = Pf \cap <H_0,...,H_5>$.
Introduce projective coordinates $(t_0:...:t_5)$
in ${\Pi}^5$ such that the point
$(t_0,...,t_5)$ represents the vector
$t_0H_0 + ... + t_5H_5$.
Then the cubic 4-fold ${\bf B}_3 = Pf \cap {\bf P}^5(t_0:...:t_5)$
is defined by
${\bf B}_3: F = 32t_0t_2t_5 -   t_0t_3t_5 -  2t_1^2t_5 - 2t_0t_4^2
+ 3t_1t_3t_4 - 12t_1t_2t_4 - 45t_2^2t_3 - 9t_2t_3^2 = 0$.

The 6-vector $(0,...,0)$
is the only value of $(t_0,...,t_5)$ where the $15$ Pfaffians
$Pf_{ij}(t_0,...,t_5)$ of the matrix
$H(t_0,...,t_5)$ = $t_0H_0 +...+t_5H_5$
vanish.
Therefore
$rank \ b = 4$ for any $b \in {\bf B}_3$.

The fourfold ${\bf B}_3$ = $(F = 0) \subset {\Pi}^5$
is \underline{singular},
and it is is not hard to check that
$Sing \ {\bf B}_3$ = $({\nabla}_{(t_0:...:t_5)}F$ =
$(0,...,0)) \subset {\Pi}^5$
is the rational normal quartic curve
$C_4$ = $\{ [r] =
(1: 2r : r^2/3 : 8r^2/3 : 2r^3 : r^4)| r \in {\bf C} \}$;
for simplicity we let $t_0 = 1$.

\medskip 

{\bf (3.18)} \
Now we are ready to prove Lemma {\bf (A)} for $g = 8$.

Let $X_{14} = G(2,6) \cap {\bf P}^9 = {\bf P}(U_{10})$
be such that $X_{14} \supset S_{14}$, and let
$B_3 = Pf \cap {\bf P}(U_{10}^{\perp})$.
Then
$S_{14} \subset X_{14} = G(2,V) \cap {\bf P}(U_{10})$
$\Leftrightarrow$
${\bf P}^8 = Span \ S_{14} \subset {\bf P}(U_{10}) = Span \ X_{14}$
$\Leftrightarrow$
${\bf P}(U_{10}^{\perp}) \subset {\Pi}^5$
$\Leftrightarrow$
$B_3 \subset {\bf B}_3$.
Since $B_3 \subset {\bf B}_3$ then
$rank(H) = 4$ for any $H \in B_3$ (see {\bf (3.17)}),
hence the Da Palatini birationalities
${\xi}: X_{14} \rightarrow B_3$ (see {\bf (3.16)})
are well-defined.

Assume that $X_{14}$ is smooth.
Then $B_3$ must be smooth (see above or [Pu, p. 83]).
But $B_3$ must be singular at any of the
intersection points of the hyperplane
${\bf P}(U_{10}^{\perp}) \subset {\Pi}^5$
and the rational quartic curve $C_4 = Sing \ {\bf B}_3$
-- contradiction (see the end of {\bf (3.17)}).
Therefore {\it any} $X_{14} \supset S_{14}$ must be singular.
{\bf q.e.d.}

\medskip 

{\bf (3.19) \ Remark.} \
If ${\bf P}(U_{10}^{\perp}) = Span(C_4)$
then $B_3$ is singular along $C_4$,
and it can be seen that then
$X_{14}$ is singular along $C_8$.
Let ${\bf P}(U_{10}^{\perp}) \not= Span(C_4)$.
Then the hyperplane ${\bf P}(U_{10}^{\perp}) \subset {\Pi}^5$
intersects the rational normal quartic $C_4 = Sing \ {\bf B}_3$
in $4$ possibly coincident points $b_1,b_2,b_3,b_4$.
Let, for simplicity $b_i$ be different from each other.
Then one can show that the general Da Palatini birationality
$B_3 \leftrightarrow X_{14}$ sends $\{ b_1,b_2,b_3,b_4 \}$ to
$4 = 12 - g(X_{14})$ singular points of $X_{14}$
which lie on $C_8$.
Let $H \subset {\bf P}^5$ be a hyperplane, and let
${\xi}^{-1}_H: B_3 \rightarrow X_{14}$ be the
Da Palatini birationality defined by $H$.
We shall see that for the general $H$ the rational map
${\xi}^{-1}_H = {\psi}_H$ is regular at a neighborhood
of any $b_i$, and ${\xi}^{-1}_H(b_i) \in C_8$.
For this, by the definition of the maps
$\phi$ and $\psi$, it is necessary to see that
the kernel-map
$ker: B_3 \rightarrow G(1:5)$, $b \mapsto n_b$
sends the quartic $C_4$ isomorphically to $C_8$.

Let $b = (b_{ij}) \in B_3$. Then the Pl\"ucker
coordinates of the line $n_b = ker(b)$ are
$(-1)^{i+j}Pf_{ij}(b)$, where $Pf_{ij}(b)$ are
the 15 quadratic Pfaffians of the skew-symmetric
matrix $\hat{b}$; note that $rank(b) = 4$ for any
$b \in B_3 \subset {\bf B}_3$.
Now, it rests only to replace
$b$ by the general point
$b(t) =
H_0 + tH_1 + t^2/12H_2 + 2t^2/3H_3 + t^3/4H_4 + t^4/16H_5 \in C_4$,
and to see that the Pl\"ucker coordinates
of $b(t)$ parameterize the general point
$x_{ij}(\vec{s})$ of $C_8$ (where $s = 2/t$) -- see {\bf (3.17)}.

\bigskip

\bigskip 

\centerline{\bf Proof of Lemma (A) for g = 7}

\medskip

{\bf (3.20)} \
In the proof of Lemma {\bf (A)} for $g = 7$ we shall need the
known by [I2] description of the projection from a line $l$
on a smooth prime Fano threefold $X_{2g-2}$ such that
$N_{l/X} = {\cal O}(1) \oplus {\cal O}(-2)$.

\medskip

{\bf (3.21) Lemma} (see \S 1 Proposition 3 in [I2]).
{\sl
Let $l$ be a line on the smooth prime Fano threefold 
$X = X_{2g-2} \subset {\bf P}^{g+1}$ such that
$N_{l/X} = {\cal O}_l(-2) \oplus {\cal O}_l(1)$, 
and let ${\sigma}:X' \rightarrow X$ be the blow up of $l$.
Let $Z' = {\sigma}^{-1}(l)$, let $H' \sim {\sigma}^*H - Z'$
be the proper preimage of the hyperplane section $H$ of $X$,
and let ${\pi}:X \rightarrow X'' \subset {\bf P}^{g-1}$
be the projection from $l$. Then:

{\bf (i).} \
If $g \ge 5$ then the composition
${\phi} = {\pi} \circ {\sigma}: X' \rightarrow {\bf P}^{g-1}$
is a birational morphism (given by the linear system $\mid H' \mid$),
to a threefold $X'' \subset {\bf P}^{g-1}$.

{\bf (ii).} \
The restriction to $Z'= {\bf P}(N_{l/V}) = {\bf F}_3$
of the linear system $\mid H' \mid$ is the complete linear system
$\mid s' + 3 f' \mid$, where $s'$ and $f'$ are the classes
of the exceptional
section and the fiber of the rational ruled surface $Z'$.

{\bf (iii).} \
The restriction ${\phi}\mid_{Z'}$ of ${\phi}$ to $Z'$ maps
$Z'$ to a cone $Z''$ over a twisted cubic curve, contracting the
exceptional section $s'$ of $Z'$ to the vertex of $Z''$.

{\bf (iv).} \
If $g \ge 7$ then there are only a finite number of lines
$l_i \subset X$ ($i = 1,...,N$) which intersect $l$.
Let $l_i' \subset X'$ of $l_i$ ($i = 1,...,N$)
be the proper preimages $l_i' \subset X'$ of $l_i$ ($i = 1,...,N$).
Then the morphism ${\phi}:X' \rightarrow X''$
is an isomorphism outside $l_1' \cup ... \cup l_N' \cup s'$,
and ${\phi}$ contracts $s'$ and any of $l_i'$ to isolated double
points of $X''$.

{\bf (v).} \
Let $H'' = {\phi}(H')$ be the hyperplane section of $X''$.
Then if $g \ge 7$ then $-K_{X''} \sim H''$, i.e. the variety
$X'' = X''_{2(g-2)-2} \subset {\bf P}^{g-1}$
is an anticanonically embedded Fano threefold with isolated
singularities as in {\bf (iv)}.
}

\medskip 

{\bf (3.22)} \
Suppose that there exists a smooth prime
$X = X_{12} \subset {\bf P}^8$
which contains the tangent scroll $S_{12}$ to the rational
normal curve $C_7$. Therefore the general such $X$ is smooth,
and we may suppose that $X \supset S_{12}$ is general.

Let $l \subset S_{12}$ be any of the tangent lines to $C_7$
and let ${\pi}:X \rightarrow X''$ be the
projection from $l$.
By [I2, \S 1] (see also {\bf (1.2)}) \
$N_{l/X} = {\cal O}(1) \oplus {\cal O}(-2)$,
therefore {\bf (3.21)(i)-(v)} take place.

By {\bf (3.21)(v)} the threefold
$X'' = X''_8 \subset {\bf P}^6$ is an anticanonically
embedded Fano threefold of genus $5$.
Let $S'' \subset X''$ be the proper image of $S_{12}$,
and let $C'' \subset S'' \subset X''$ be the proper image
of $C_7$. It is easy to see that $C'' = C''_5$ is a rational
normal quintic, and $S'' = S''_8 \subset X''$ is the tangent
scroll to $C''$.

In order to use the proof of Lemma {\bf (A)} for $g = 5$
we have to see whether $X''_8 \subset {\bf P}^6$
is, in fact, a complete intersection of three quadrics.
If $X''$ were nonsingular then the classification
of the smooth Fano threefolds will imply that $X''$
will be a complete intersection of three quadrics.
But $X''$ is singular -- see {\bf (3.21)(iv),(v)}.

However, especially in this case, $X''$ = $X''_8 \subset {\bf P}^6$
is {\it still} a complete intersection of three quadrics
(see Theorem (6.1) (vii) in [I1]).

Denote by $Sing(X)\cap{C}$ the set of singular points
of $X = X_{12}$ on $C$,
and let $Sing(X'')\cap{C''}$ be the set of singular
points of $X''$ on $C''$.

\medskip 

{\bf (3.23)} \
By the proof of Lemma {\bf (A)} for $g = 5$,
the elements of $Sing(X'')\cap{C}$ are in a
(1:1) correspondence
with the different zeros of a homogeneous form
$F_7(s_o:s_1)$ of degree $7$ (see {\bf (3.6)}).
Clearly, the vertex $o$ of $Z''$ lies on $C''$.
Moreover, by {\bf (3.21)(iv)}, $o$ is a double singularity of $X''$.
Since $l \cap Sing \ X$ = $\emptyset$, and since the
tangent lines $l' \not= l$ to $C = C_7$ do not intersect $l$,
{\bf (3.21)(iv)} yield that, set-theoretically:
$Sing(X)\cap{C}$ $\cong$ $Sing(X'')\cap {C''} - \{ o \}$.

Let $F_7(s_0:s_1) = 0$ be as above.
Since $dim \ Sing \ X'' = 0$,
the form $F_7$ does not vanish
on $C_7$; and we can assume that o = (1:0) and
$F_7(0:1) \not= 0$. Therefore if
$s = s_1/s_0$ and $f_7(s) = F(1:s)$ then
$deg \ f_7(s) = 7$.
By the previous the elements of
$Sing(X)\cap{C}$ correspond to the different zeros
of the polynomial
$s^{-m}.f_7(s) = 0$, where $m = mult_o \ f_7(s)$.

By the local definition of $m$,
the integer $m = m(o) = mult_o f_7(s)$
does not depend on the genus $g \ge 7$
of $X_{2g-2}$
as well on the choice of the general
tangent line $l$ to $C_g$.
It can be seen that $m = 2$,
but for the proof it is enough
to know that $m \le 2$.

\medskip

{\bf $({\bf \ast})$ Lemma}. \ $m \le 2$.

\medskip

{\bf Proof of} $({\ast}).$ \
By construction $X'' \supset S'' \cup Z''$
where
$S'' = S_8''$ is the tangent scroll to the rational normal
quintic $C'' = C_5'$ such that $o \in C''$,
$Z''$ is a cone over a twisted cubic, and $o$ is the vertex of $Z$.
Moreover $Z''$ is triple tangent
to $S''$ at the tangent line $F$ to $C''$ at $o$.
Indeed $S''$ is a hyperplane section of $Z''$
which passes through the vertex $o$ of $Z$.
Therefore $S''.Z = f_1+f_2+f_3$ is a sum of 3 rulings of $Z''$.
Since $f_i$ are rulings of $Z$, $o \in f_i$ for $i = 1,2,3$.
Therefore any $f_i$ is a line on the tangent scroll $S''$ to $C''$
which passes through $o \in C''$. Therefore $f_i = F$ must be
a tangent line to $C''$ at $o$, i.e. $S''.Z'' = 3F$.

By Theorem 9.9 in [I3], the general
complete intersection $X_8 \subset {\bf P}^6$ of three quadrics,
containing a cone $Z_3$ over a twisted cubic, is a projection
of $X_{12}$ from a line $l$ such that
$N_{l/X_{12}} = {\cal O}(1) \oplus {\cal O}(-2)$.
The inverse of the projection ${\pi}_l$ is defined by
the linear system $\mid H + Z_3 \mid$, where $H$ is the hyperplane
section of $X_8$.

Let $X_8 \supset Z_3 \cup S_8$ be as above.
Then $X_{12}$ will contain a tangent scroll $S_{12}$
to a rational normal curve $C_7$, and $l$ will be a tangent
line to $C_7$.
Therefore any $X_8 \supset Z_3 \cup S_8$
will be a deformation of a projection
of $X_{12} \subset S_{12}$
from a tangent line to $C_7$.

It rests to see that $m(X_8) = mult_o f_7 \le 2$
for $f_7$ corresponding, as above, to some particular
such $X_8$.

\medskip

{\sc Example.} \
Let ${\bf P}^5(x) = {\bf P}^5(x_0:...:x_5)$, and let
$q_0 = -x_0x_4 + 4x_1x_3 - 3x_2^2$, \
$q_1 = -x_0x_5 + 3x_1x_4 - 2x_2x_3$, \
$q_2 = -x_1x_5 + 4x_2x_4 - 3x_3^2$.
Then
$S_8 = (q_0 = q_1 = q_2 = 0) \subset {\bf P}^5(x)$
will be the tangent scroll to the
rational normal quintic
$C_5: x_i = s_0^{5-i}s_1^i (0 \le i \le 5)$.

Let
$X_8 = (Q_0 = Q_1 = Q_2 = 0)
\subset {\bf P}^6(x:u) = {\bf P}^6(x_0:x_1:x_2:x_3:x_4:x_5:u)$,
where

\medskip

\hspace{3cm}{$Q_0 = q_0 + L_o(x_4:x_5)u$,}

\hspace{3cm}{$Q_1 = q_1 + (12x_1 + L_1(x_4:x+5))u$,}

\hspace{3cm}{$Q_2 = q_2 + (27/2 x_2 + L_2(x_4:x_5))u$,}

\medskip
\noindent
$L_0$, $L_1$ and $L_2$ being linear forms of $(x_4:x_5)$.
Evidently $X_8 \cap (u = 0) = S_8$.

Let ${\bf P}^4 = {\bf P}^4(x_0:x_1:x_2:x_3:u) \subset {\bf P}^6$,
and let $Z_3 = X_8 \cap {\bf P^4}$.
Then $Z_3 = (P_0 = P_1 = P_2 = 0) \subset {\bf P}^4$, where
$P_0 = x_1x_3/3 - x_2^2/4$, \
$P_1 = x_1u - x_2x_3/6$, \
$P_2 = x_2u - x_3^2/9$.

Therefore $Z_3$ is a cone with center $o = (1:0:...:0) \in C_5$
over the twisted cubic curve
$C_3 = Z_3 \cap (x_0 = 0)$,
$C_3: (x_1:x_2:x_3:u) = (t_0^3:2t_0^2t_1:3t_0t_1^2:t_1^3)$.
Let $s = s_1/s_0$, and we may suppose that the
point $(0:...:0:1) \in C_5$ is not a singular point of $X_8$.
Then, by {\bf (3.6)}, the equation of $(Sing \ X_8)|_{C_5}$ is

\medskip

${}$\hspace{1cm}{$f_7(s) =
s^2  \partial Q_0 / \partial u (1:s:...:s^5)
- s \partial Q_1 / \partial u (1:s:...:s^5)
+ \partial Q_2 / \partial u (1:s:...:s^5)$}

${}$\hspace{1cm}{= $s^2L_0(s^4,s^5) - s(12s + L_1(s^4,s^5))
+ (27/2s^2 + L_2(s^4,s^4)) = 3/2s^2 + o(s^3)$,}

\medskip
\noindent
where $o(s^3)$ is a sum of terms of degree $\ge 3$.
Therefore $m(X_8) = mult_o f_7(s) = 2$. {\bf q.e.d.}

\medskip 

{\bf (3.24)} \
Let $X$ = $X_{12} \supset S_{12}$ be general.
Since $m \le 2$ then
$deg \ s^{-m}f_7(s)$
$\ge 7 - 2 = 5$ = $12 - g(X_{12}) > 0$.
In particular $g(s) := s^{-m}f_7(s)$ is not a constant.
Since $g(0) \not= 0$,
and since the elements of $Sing(X) \cap C$ are
in a (1:1) correspondence with the different
zeros of $g(s) = s^{-m}f_7(s)$ (see above),
then $X_{12}$ must be singular,
which contradicts the initial assumption.
This proves Lemma {\bf (A)} for $g = 7$.

\bigskip

\centerline{\bf Proof of Lemma (A) for g = 9.}

\medskip

{\bf (3.25)} \
Let $X_{16} \subset {\bf P}^{10}$ contains the tangent scroll
$S = S_{16}$ to the rational normal curve $C = C_9$, and
suppose that nevertheless $X_{16}$ is smooth.

\medskip

Let $L \subset X_{16}$ be a tangent line to $C$,
and consider the double projection ${\pi} = {\pi}_{2L}$ of $X$
from the line $L$, i.e. ${\pi}$ is the rational map on $X$ defined
by the non-complete linear system $|{\cal O}_X(1 - 2L)|$.
Since $X = X_{16}$ is assumed to be smooth then,
by \S 2 in [I2]:

\medskip

${\bf ({\ast})}$. \ 
{\sl 
${\pi} = {\pi}_{2L}$ sends $X$ birationally to ${\bf P}^3$.
Moreover, on ${\bf P}^3$ there exists
a smooth irreducible curve $C = C^3_7$ of genus $3$ and degree $7$,
which lies on a unique cubic surface $S_3 \subset {\bf P}^3$, 
and such that the inverse to ${\pi}$ birational map
${\phi}: {\bf P}^3 \rightarrow X$ is given by the non-complete
linear system $|{\cal O}_{{\bf P}^3}(7 - 2C)|$. 
} 

\medskip 

By ${\bf ({\ast})}$, the proper image ${\pi}(H)$
of any hyperplane section $H \subset X$
is an irreducible component of an effective divisor 
$S_7 \in |{\cal O}_{{\bf P}^3}(7 - 2C)|$.
If moreover $H$ contains the line $L$
but $H \not\in |{\cal O}_X(1 - 2L)|$
(for example if $H = S_{16}$) 
then ${\pi}(H) \subset {\bf P}^3$ will be a quartic surface
containing the curve $C = C^3_7$ (see the proof of the
Main Theorem in \S 2 of [I2]), and in this case
$S_7 = {\pi}(H) + S_3 \in |{\cal O}_{{\bf P}^3}(7 - 2C)|$. 

\smallskip 

Therefore $S_4 := {\pi}(S_{16})$ is a quartic surface
in ${\bf P}^3$ containing the curve $C = C^3_7$.
Moreover, the double projection
$\pi$ sends the general tangent line
$L'$ to $C_9$ to a tangent line
${\pi}(L')$ to the proper image ${\pi}(C_9)$;
and since $C_9 \cong {\bf P}^1$ then ${\pi}(C_9)$
is rational. 
Therefore the quartic surface $S_4 = {\pi}(S_{16})$
is the tangent scroll to the rational curve
${\pi}(C_9) \subset{\bf P}^3$. 
The last is only possible if ${\pi}(C_9) = C_3$ is a
twisted cubic and $S_4$ is the tangent scroll to $C_3$, 
and we shall see that this is impossible.

\smallskip

The surface $S_4 \subset {\bf P}^3$
is the tangent scroll to the twisted cubic $C_3$.
Then, by Lemma 1.6 and p. 498 in [MU],
the normalization of $S_4$
is the quadric ${\bf P}^1 \times {\bf P}^1$,
and the map ${\nu}: {\bf P}^1 \times {\bf P}^1 \rightarrow S_4$
is given by a linear system of bidegree $(1,2)$. 

Let $\Gamma \subset {\bf P}^1 \times {\bf P}^1$
be the proper transform of $C^3_7$, and let $(a,b)$ be the
bidegree of $\Gamma$. Therefore $7 = deg(C^3_7) = 2a + b$,
and $3 = g(C^3_7) = g({\Gamma}) = (a-1)(b-1)$.
Obviously, these two equations for the integers $a$ and $b$
have no integral solutions -- contradiction.

Therefore $X_{16} \supset S_{16}$ can't be smooth,
which proves Lemma {\bf (A)} in case $g = 9$.

\bigskip

This completes the proof of Lemma {\bf (A)}.

\medskip
\noindent
{\it Carmen Schuhmann} \ \
University of Leiden, \ P.O. Box 9512\\
${}$ \hspace{3.6cm} 2300 RA Leiden, \
The Netherlands

\medskip
\noindent
{\it Atanas Iliev} \ \
Institute of Mathematics, Bulgarian Academy of Sciences\\
${}$ \hspace{2.2cm} Acad. G. Bonchev Str. 8, \ 1113 Sofia, \
Bulgaria

\end{document}